\documentclass{conm-p-l}  
\usepackage{amssymb}

\newtheorem{theorem}{Theorem}
\newtheorem{proposition}[theorem]{Proposition}
\newtheorem{definition}[theorem]{Definition}
\newtheorem{lemma}[theorem]{Lemma}
\newtheorem{corollary}[theorem]{Corollary}

\numberwithin{equation}{section}
\numberwithin{theorem}{section}

\begin{document}
\title[Relative tensor products]{Relative Tensor Products for Modules over von Neumann Algebras}
\author{David Sherman}
\address{Department of Mathematics\\ University of Illinois\\ Urbana, IL 61801-2975}
\email{dasherma@math.uiuc.edu}

\begin{abstract}

We give an overview of relative tensor products (RTPs) for von Neumann algebra modules.  For background, we start with the categorical definition and go on to examine its algebraic formulation, which is applied to Morita equivalence and index.  Then we consider the analytic construction, with particular emphasis on explaining why the RTP is not generally defined for every pair of vectors.  We also look at recent work justifying a representation of RTPs as composition of unbounded operators, noting that these ideas work equally well for $L^p$ modules.  Finally, we prove some new results characterizing preclosedness of the map $(\xi, \eta) \mapsto \xi \otimes_\varphi \eta.$

\end{abstract}

\maketitle

\section{Introduction}

The purpose of this article is to summarize and explore some of the various constructions of the \textit{relative tensor product} (RTP) of von Neumann algebra modules.  Alternately known as composition or fusion, RTPs are a key tool in subfactor theory and the study of Morita equivalence.  The idea is this: given a von Neumann algebra $\mathcal{M}$, we want a map which associates a vector space to certain pairs of a right $\mathcal{M}$-module and a left $\mathcal{M}$-module.  If we write module actions with subscripts, we have
$$(\mathfrak{X}_\mathcal{M}, {}_\mathcal{M} \mathfrak{Y}) \mapsto \mathfrak{X} \otimes_\mathcal{M} \mathfrak{Y}.$$
This should be functorial, covariant in both variables, and appropriately normalized.  Other than this, we only need to specify which modules and spaces we are considering.

In spirit, RTPs are algebraic; a ring-theoretic definition can be found in most algebra textbooks.  But in the context of operator algebras, the requirement that the output be a certain type of space - typically a Hilbert space - causes an analytic obstruction.  As a consequence, there are domain issues in any vector-based construction.  Fortunately, von Neumann algebras have a sufficiently simple representation theory to allow a recasting of RTPs in algebraic terms.

The analytic study of RTPs can be related nicely to noncommutative $L^p$ spaces.  Indeed, examination of the usual ($L^2$) case reveals that the technical difficulties come from a ``change of density".  (We say that the \textit{density} of an $L^p$-type space is $1/p$.)  Once this is understood, it is easy to handle $L^p$ modules [JS] as well.  Modular algebras ([Y], [S]) provide an elegant framework, so we briefly explain their meaning.

The final section of the paper investigates the question, ``When is the map $(\xi, \eta) \mapsto \xi \otimes_\varphi \eta$ preclosed?"  This may be considered as an extension of Falcone's theorem [F], in which he found conditions for the map to be everywhere-defined.  We consider a variety of formulations.

We have tried to make the paper as accessible as possible to non-operator algebraists, especially in the first half.  Of course, even at this level many results rely on familiarity with the projection theory of von Neumann algebras; basic sources are [T1], [T2], [KR].  Primary references for RTPs are [Sa], [P], [F], [C2].

\section{Notations and background}

The basic objects of this paper are \textit{von Neumann algebras}, always denoted here by $\mathcal{M},$ $\mathcal{N},$ or $\mathcal{P}.$  These can be defined in many equivalent ways:
\begin{itemize}
\item C*-algebras which are dual spaces.
\item strongly-closed unital *-subalgebras of $\mathcal{B}(\mathfrak{H})$.  $\mathcal{B}(\mathfrak{H})$ is the set of bounded linear operators on a Hilbert space $\mathfrak{H}$; the strong topology is generated by the seminorms $x \mapsto \|x \xi \|, \xi \in \mathfrak{H}$; the * operation is given by the operator adjoint.
\item *-closed subsets of $\mathcal{B}(\mathfrak{H})$ which equal their double (iterated) commutant.  The commutant of a set $S \subset \mathcal{B}(\mathfrak{H})$ is $\{x \in \mathcal{B}(\mathfrak{H}) \mid xy = yx, \forall y \in S\}$.
\end{itemize}
As one might guess from the definitions, the study of von Neumann algebras turns on the interplay between algebraic and analytic techniques.

Finite-dimensional von Neumann algebras are direct sums of full matrix algebras.  At the other extreme, commutative von Neumann algebras are all of the form $L^\infty(X, \mu)$ for some measure space $(X, \mu)$, so the study of general von Neumann algebras is considered ``noncommutative measure theory."  Based on this analogy, the (unique) predual $\mathcal{M}_*$ of $\mathcal{M}$ is called $L^1(\mathcal{M})$; it is the set of normal (= continuous in yet another topology, the $\sigma$-weak) linear functionals on $\mathcal{M} \subset \mathcal{B}(\mathfrak{H})$, and can be thought of as ``noncommutative countably additive measures".  A functional $\varphi$ is positive when $x > 0 \Rightarrow \varphi(x) \ge 0$; the set of positive normal functionals is denoted $\mathcal{M}_*^+.$  The \textit{support} $s(\varphi)$ of a positive normal linear functional $\varphi$ is the smallest projection $q \in \mathcal{M}$ with $\varphi(1-q) = 0.$  So if $\mathcal{M}$ is abelian, $\varphi$ corresponds to a measure and $q$ is the (indicator function of the) usual support.

For simplicity, all modules in this paper are separable Hilbert spaces (except in Section \ref{S:modalg}), all algebras have separable predual, all linear functionals are normal, and all representations are normal and nondegenerate ($\mathcal{M}\mathfrak{H}$ or $\mathfrak{H}\mathcal{M}$ is all of $\mathfrak{H}$).  Two projections $p,q$ in a von Neumann algebra are said to be (Murray-von Neumann) \textit{equivalent} if there exists $v \in \mathcal{M}$ with $v^*v =p$, $vv^* = q$.  Such an element $v$ is called a partial isometry, and we think of $p$ and $q$ as being ``the same size".  Subscripts are used to represent actions, so $\mathfrak{X}_\mathcal{M}$ indicates that $\mathfrak{X}$ is a right $\mathcal{M}$-module, i.e. a representation of the opposite algebra $\mathcal{M}^{\text{op}}.$  It is implicit in the term ``bimodule", or in the notation ${}_\mathcal{M} \mathfrak{H}_\mathcal{N}$, that the two actions commute.  The phrase ``left (resp. right) action of" is frequently abbreviated to $L$ (resp. $R$) for operators or entire algebras, so that we speak of $L(x)$ or $R(\mathcal{M})$.  Finally, we often write $M_\infty$ for the von Neumann algebra of all bounded operators on a separable infinite-dimensional Hilbert space, and $M_\infty(\mathcal{M})$ for the von Neumann tensor product $M_\infty \bar{\otimes} \mathcal{M}$.  One can think of this as the set of infinite matrices with entries in $\mathcal{M}$; we will denote by $e_{ij}$ the matrix unit with 1 in the $ij$ position and 0 elsewhere.

The (left) representation theory of von Neumann algebras on Hilbert spaces is simple, so we recall it briefly.  (Most of this development can be found in Chapters 1 and 2 of [JoS].)  First, there is a standard construction, due to Gelfand-Neumark and Segal (abbreviated GNS), for building a representation from $\varphi \in \mathcal{M}_*^+$.  To each $x \in \mathcal{M}$ we formally associate the vector $x\varphi^{1/2}$ (various notations are in use, e.g. $\eta_\varphi(x)$ or $\Lambda_\varphi(x)$, but this one is especially appropriate ([C2] V.App.B, [S])).  We endow this set with the inner product
$$<x \varphi^{1/2} , y\varphi^{1/2}> = \varphi(y^* x),$$
and set $\mathfrak{H}_\varphi$ to be the closure in the inherited topology, modulo the null space.  The left action of $\mathcal{M}$ on $\mathfrak{H}_\varphi = \overline{\mathcal{M} \varphi^{1/2} }$ is bounded and densely defined by left composition.  

When $\varphi$ is faithful (meaning $x>0 \Rightarrow \varphi(x)>0$), the vector $\varphi^{1/2} = 1\varphi^{1/2}$ is cyclic ($\overline{\mathcal{M} \varphi^{1/2}} = \mathfrak{H}_\varphi$) and separating ($x \ne 0 \Rightarrow x\varphi^{1/2} \ne 0$).  Now all representations with a cyclic and separating vector are isomorphic - a sort of ``left regular representation"; we will denote this by ${}_\mathcal{M} L^2(\mathcal{M}).$  It is a fundamental fact that the commutant of this action is antiisomorphic to $\mathcal{M}$, and when we make this identification we call ${}_\mathcal{M} L^2(\mathcal{M})_\mathcal{M}$ the \textit{standard form} of $\mathcal{M}$.  If $\varphi$ is not faithful, the GNS construction produces a vector $\varphi^{1/2}$ which is cyclic but not separating, and a representation which is isomorphic to ${}_\mathcal{M} L^2(\mathcal{M}) s(\varphi)$ ([T2], Ch. VIII, IX).  

Now let us examine an arbitrary (separable, so countably generated) module ${}_\mathcal{M}\mathfrak{H}$.  Following standard arguments (e.g. [T1] I.9), $\mathfrak{H}$ decomposes into a direct sum of cyclic representations ${}_\mathcal{M}(\overline{\mathcal{M}\xi_n }) $, each of which is isomorphic to the GNS representation for the associated vector state $\omega_{\xi_n}$ ($=<\cdot \:\xi_n, \xi_n>$).  With $q_n = s(\omega_{\xi_n})$, we have
$${}_\mathcal{M} \mathfrak{H} \simeq \bigoplus {}_\mathcal{M} \overline{\mathcal{M} \xi_n } \simeq \bigoplus {}_\mathcal{M} \mathfrak{H}_{\omega_{\xi_n}} \simeq \bigoplus {}_\mathcal{M} 
L^2(\mathcal{M}) q_n.$$
(Here and elsewhere, ``$\simeq$" means a unitary equivalence of (bi)modules.)  Since this is a left module, it is natural to write vectors as rows with the $n$th entry in $L^2(\mathcal{M}) q_n$:
\begin{equation} \label{E:rowrep}
\mathfrak{H} \simeq (L^2(\mathcal{M})q_1 \: L^2(\mathcal{M}) q_2 \: \cdots ) \simeq (L^2(\mathcal{M}) \: L^2(\mathcal{M}) \: \cdots) (\sum q_n \otimes e_{nn}).
\end{equation}
We will call such a decomposition a \textit{row representation} of ${}_\mathcal{M}\mathfrak{H}$.  Here $e_{nn}$ are diagonal matrix units in $M_\infty$, so $(\sum q_n \otimes e_{nn})$ is a diagonal projection in $M_\infty(\mathcal{M}).$  The left action of $\mathcal{M}$ is, of course, matrix multiplication (by $1\times 1$ matrices) on the left.

The module $(L^2(\mathcal{M}) \: L^2(\mathcal{M}) \: \cdots)$ will be denoted $R^2(\mathcal{M})$ (for ``row").  Since the standard form behaves naturally with respect to restriction - $L^2(q\mathcal{N}q) \simeq qL^2(\mathcal{N})q$ as bimodules - it follows that $L^2(M_\infty(\mathcal{M}))$ is built as infinite matrices over $L^2(\mathcal{M})$ (see \eqref{E:l2l2}).  Thus $R^2(\mathcal{M})$ can be realized as $e_{11} L^2(M_\infty(\mathcal{M})).$

\begin{proposition} \label{T:module}
Any countably generated left representation of $\mathcal{M}$ on a Hilbert space is isomorphic to $R^2(\mathcal{M})q$ for some diagonal projection $q \in M_\infty(\mathcal{M}).$  Any projection in $M_\infty(\mathcal{M})$, diagonal or not, defines a module in this way, and two such modules are isomorphic exactly when the projections are equivalent.  In fact
\begin{equation} \label{E:comm}
\text{Hom} ({}_\mathcal{M} R^2(\mathcal{M})q_1, {}_\mathcal{M} R^2(\mathcal{M}) q_2) = R(q_1 M_\infty(\mathcal{M})q_2).
\end{equation}
So isomorphism classes correspond to equivalence classes of projections in $M_\infty(\mathcal{M})$, which is the monoid $V(M_\infty(\mathcal{M}))$ in $K$-theoretic language [W-O].  The direct sum of isomorphism classes of modules corresponds to the sum of orthogonal representatives in $V(M_\infty(\mathcal{M}))$, giving a monoidal equivalence.
\end{proposition}

We denote the category of separable left $\mathcal{M}$-modules by $Left\: L^2(\mathcal{M}).$  For us, the most important consequence of \eqref{E:comm} is that
\begin{equation} \label{E:rcomm}
\mathcal{L}({}_\mathcal{M} R^2(\mathcal{M}) q) = R(qM_\infty(\mathcal{M})q),
\end{equation}
where ``$\mathcal{L}$" stands for the commutant of the $\mathcal{M}$-action.  (In particular, the case $q=e_{11}$ is just the standard form.)  The algebra $qM_\infty(\mathcal{M})q$ is called an \textit{amplification} of $\mathcal{M}$, being a generalization of a matrix algebra with entries in $\mathcal{M}$.  
Of course everything above can be done for right modules - the relevant abbreviations are $C^2(\mathcal{M}),$ for ``column," and $Right\:L^2(\mathcal{M})$.

\textbf{Example.} Suppose $\mathcal{M}=M_3(\mathbb{C})$.  In this case the standard form may be taken as
$${}_{M_3}L^2(M_3)_{M_3}; \quad L^2(M_3) \simeq (M_3, <\cdot, \cdot>), \text{ where }<x,y> = \text{Tr}(y^*x).$$
Note that this norm, called the Hilbert-Schmidt norm, is just the $\ell^2$ norm of the matrix entries, and that the left and right multiplicative actions are commutants.  (If we had chosen a nontracial positive linear functional, we would have obtained an isomorphic bimodule with a ``twisted" right action... this is inchoate Tomita-Takesaki theory.)  The module $R^2(M_3)$ is $M_{3 \times\infty}$, again with the Hilbert-Schmidt norm, and the commutant is $M_\infty(M_3) \simeq M_\infty.$  According to Proposition \ref{T:module}, isomorphism classes of left $M_3$-modules should be parameterized by equivalence classes of projections in $M_\infty$.  These are indexed by their rank $n \in (\mathbb{Z}_+ \cup \infty)$; the corresponding isomorphism class of modules has representative $M_{3\times n}$.  In summary, we have learned that any left representation of $M_3$ on a Hilbert space is isomorphic to some number of copies of $\mathbb{C}^3$.  The same argument shows that $V(M_\infty(M_k)) \simeq (\mathbb{Z}_+ \cup \infty)$ for any $k$.

Properties of the monoid $V(M_\infty(\mathcal{M}))$ determine the so-called \textit{type} of the algebra.  For a \textit{factor} (a von Neumann algebra whose center is just the scalars), there are only three possibilities: $(\mathbb{Z}_+ \cup \infty),$ $(\mathbb{R}_+ \cup \infty),$ and $\{0, +\infty\}$.  These are called types I, II, III, respectively; a fuller discussion is given in Section \ref{S:preclosed}.

\section{Algebraic approaches to RTPs} \label{S:algebraic}

When $R$ is a ring, the algebraic $R$-relative tensor product is the functor, covariant in both variables, which maps a right $R$-module $A$ and left $R$-module $B$ to the vector space $(A \otimes_{alg} B) / N$, where $N$ is the subspace generated algebraically by tensors of the form $ar \otimes b - a \otimes rb$.  In functional analysis, where spaces are usually normed and infinite-dimensional, one obvious amendment is to replace vector spaces with their closures.  But in the context of Hilbert modules over a von Neumann algebra $\mathcal{M}$, this is still not enough.  Surprisingly, a result of Falcone ([F], Theorem 3.8) shows that if the RTP $L^2(\mathcal{M}) \otimes_\mathcal{M} L^2(\mathcal{M})$ is the closure of a continuous (meaning $\|I(\xi \otimes \eta)\| < C \|\xi\| \|\eta\|$) nondegenerate image of the algebraic $\mathcal{M}$-relative tensor product, $\mathcal{M}$ must be \textit{atomic}, i.e. $\mathcal{M} \simeq \oplus_n \mathcal{B}(\mathfrak{H}_n).$  We will discuss the analytic obstruction further in Section \ref{S:analytic}.  For now, we take Falcone's theorem as a directive: do not look for a map which is defined for every pair of vectors.  If we give up completely on a vector-level construction, we can at least make the functorial

\begin{definition}
[Sa] Given a von Neumann algebra $\mathcal{M}$, a relative tensor product is a functor, covariant in both variables,
\begin{equation}
Right\, L^2(\mathcal{M}) \times Left\, L^2(\mathcal{M}) \to Hilbert: \quad (\mathfrak{H}, \mathfrak{K}) \mapsto \mathfrak{H} \otimes_{\mathcal{M}} \mathfrak{K},
\end{equation}
which satisfies
\begin{equation} \label{E:normalized}
L^2(\mathcal{M}) \otimes_{\mathcal{M}} L^2(\mathcal{M}) \simeq L^2(\mathcal{M})
\end{equation}
as bimodules. 
\end{definition}

Although at first glance this definition seems broad, in fact we see in the next proposition that there is exactly one RTP functor (up to equivalence) for each algebra.  The reader is reminded that functoriality implies a mapping of intertwiner spaces as well, so it is enough to specify the map on representatives of each isomorphism class.  In particular we have the bimodule structure
$${}_{\mathcal{L}(\mathfrak{H}_{\mathcal{M}})} (\mathfrak{H} \otimes_{\mathcal{M}} \mathfrak{K})_{\mathcal{L}({}_{\mathcal{M}}\mathfrak{K})}.$$

\begin{proposition} \label{T:colxrow}
Let $\mathfrak{H} \simeq p\: C^2(\mathcal{M}) \in Right \:L^2 Mod(\mathcal{M})$ and $\mathfrak{K} \simeq  R^2(\mathcal{M}) q \in$ $ Left\:L^2 Mod(\mathcal{M})$
for some projections $p,q \in M_\infty(\mathcal{M}).$  Then 
$$\mathfrak{H} \otimes_\mathcal{M} \mathfrak{K} \simeq 
p\:L^2(M_\infty(\mathcal{M}))q$$
with natural action of the commutants.
\end{proposition}

\begin{proof}
By implementing an isomorphism, we may assume that the projections are diagonal: $p = \sum p_i \:\otimes e_{ii}, \: q = \sum q_j \otimes e_{jj}.$  Using \eqref{E:normalized} and functoriality, we have the bimodule isomorphisms
$$\mathfrak{H} \otimes_\mathcal{M} \mathfrak{K} \simeq 
\left(\oplus p_i \:L^2(\mathcal{M}) \right) \otimes_\mathcal{M} \left( \oplus L^2(\mathcal{M}) q_j \right)$$
$$ \simeq \bigoplus_{i,j} p_i \: L^2(\mathcal{M}) \otimes_\mathcal{M} L^2(\mathcal{M})q_j \simeq \bigoplus_{i,j} p_i \:L^2(\mathcal{M}) q_j \simeq p \: L^2(M_\infty(\mathcal{M})) q.$$
\end{proof}

Visually,
\begin{equation} \label{E:l2l2}
\left( (p)
\left(
\begin{smallmatrix}
L^2(\mathcal{M})\\
L^2(\mathcal{M})\\
\vdots
\end{smallmatrix}
\right),\:
\left(
\begin{smallmatrix}
L^2(\mathcal{M}) & L^2(\mathcal{M}) & \dots
\end{smallmatrix}
\right)
(q) \right) 
\mapsto (p) \left(
\begin{smallmatrix}
L^2(\mathcal{M}) & L^2(\mathcal{M}) & \dots\\
L^2(\mathcal{M}) & L^2(\mathcal{M}) & \dots\\
\vdots & \vdots & \ddots\\
\end{smallmatrix}
\right) (q),
\end{equation}
where of course the $\ell^2$ sums of the norms of the entries in these matrices are finite.

\bigskip

After making the categorical definition above, Sauvageot immediately noted that it gives us no way to perform computations.  We will turn to his analytic construction in Section \ref{S:analytic}; here we discuss an approach to bimodules and RTPs due to Connes.  In his terminology a bimodule is called a \textit{correspondence}.  (The best references known to the author are [C2] and [P], but there was an earlier unpublished manuscript which is truly the source of Connes fusion.)

Consider a correspondence ${}_\mathcal{M} \mathfrak{H}_\mathcal{N}.$  Choosing a row representation $R^2(\mathcal{M})q$ for $\mathfrak{H}$, we know that the full commutant of $L(\mathcal{M})$ is isomorphic to $R(q M_\infty(\mathcal{M}) q)$. This gives us a unital injective *-homomorphism $\rho: \mathcal{N} \hookrightarrow q M_\infty(\mathcal{M}) q$, and from the map $\rho$ one can reconstruct the original bimodule (up to isomorphism) as ${}_\mathcal{M} (R^2(\mathcal{M})q)_{\rho(\mathcal{N})}$.

What if we had chosen a different row representation $R^2(\mathcal{M})q'$ and obtained $\rho': \mathcal{N} \to q'M_\infty(\mathcal{M})q'$?  By Proposition \ref{T:module}, the module isomorphism
\begin{equation} \label{E:rightiso}
{}_\mathcal{M} R^2(\mathcal{M})q \simeq {}_\mathcal{M} R^2(\mathcal{M})q'
\end{equation}
is necessarily given by the right action of a partial isometry $v$ between $q$ and $q'$ in $M_\infty(\mathcal{M})$.  Then $\rho$ and $\rho'$ differ by an inner perturbation: $\rho(x) = v^*\rho'(x)v$.  We conclude that \textit{the class of $\mathcal{M} - \mathcal{N}$ correspondences, modulo isomorphism, is equivalent to the class of unital injective *-homomorphisms from $\mathcal{N}$ into an amplification of $\mathcal{M}$, modulo inner perturbation.}  (These last are called \textit{sectors} in subfactor theory.)  Our convention here is to use the term ``correspondence" to mean a representative *-homomorphism for the bimodule.  But the reader should be warned that the distinction between bimodules, morphisms, and their appropriate equivalence classes is frequently blurred in the literature, sometimes misleadingly.  

Notice that a unital inclusion $\mathcal{N} \subset \mathcal{M}$ gives the bimodule ${}_\mathcal{M} L^2(\mathcal{M})_\mathcal{N}.$

\bigskip

The RTP of correspondences is extremely simple.

\begin{proposition} \label{T:corresp}
Consider bimodules ${}_\mathcal{M} \mathfrak{H}_\mathcal{N}$ and ${}_\mathcal{N} \mathfrak{K}_\mathcal{P}$ coming from correspondences $\rho_1: \mathcal{N} \hookrightarrow q M_\infty(\mathcal{M}) q$ and $\rho_2: \mathcal{P} \hookrightarrow q' M_\infty(\mathcal{N}) q'$.  The bimodule ${}_\mathcal{M}(\mathfrak{H} \otimes_\mathcal{N} \mathfrak{K})_\mathcal{P}$ is the correspondence $\rho_1 \circ \rho_2$, where we amplify $\rho_1$ appropriately.
\end{proposition}

We pause to mention that it is also fruitful to realize correspondences in terms of completely positive maps.  We shall have nothing to say about this approach; the reader is referred to [P] for basics or [A2] for a recent investigation.

\section{Applications to Morita equivalence and index}

An invertible *-functor from $Left\:L^2 Mod (\mathcal{N})$ to $Left\:L^2 Mod(\mathcal{M})$ is called a \textit{Morita equivalence} [R].  Here a \textit{*-functor} is a functor which commutes with the adjoint operation at the level of morphisms.  One way to create *-functors is the following: to the bimodule ${}_\mathcal{M} \mathfrak{H}_\mathcal{N}$, we associate
\begin{equation} \label{E:rtpfunctor}
F_\mathfrak{H}: Left\:L^2 Mod (\mathcal{N}) \to Left\:L^2 Mod(\mathcal{M}); \qquad {}_\mathcal{N} 
\mathfrak{K} \longmapsto \left( {}_\mathcal{M} \mathfrak{H}_\mathcal{N} \right) 
\otimes_\mathcal{N} \left({}_\mathcal{N} \mathfrak{K} \right).
\end{equation}

The next theorem is fundamental.

\begin{theorem} \label{T:morita}

When $L(\mathcal{M})$ and $R(\mathcal{N})$ are commutants on $\mathfrak{H}$, the RTP functor $F_\mathfrak{H}$ is a Morita equivalence.  Moreover, every Morita equivalence is equivalent to an RTP functor.

\end{theorem}

This type of result - the second statement is an operator algebraic analogue of the Eilenberg-Watts theorem - goes back to several sources, primarily the fundamental paper of Rieffel [R].  His investigation was more general and algebraic, and his bimodules were not Hilbert spaces but rigged self-dual Hilbert C*-modules, following Paschke [Pa].  From a correspondence point of view, rigged self-dual Hilbert C*-modules and Hilbert space bimodules give the same theory; the equivalence is discussed nicely in [A1].  (And the former is nothing but an $L^\infty$ version of the latter, as explained in [JS].)  Our Hilbert space approach here is parallel to that of Sauvageot [Sa], and streamlined by our standing assumption of separable preduals.

We will need 

\begin{definition}
The \textbf{contragredient} of the bimodule ${}_\mathcal{M} 
\mathfrak{H}_\mathcal{N}$ is the bimodule ${}_\mathcal{N} 
\bar{\mathfrak{H}}_\mathcal{M}$, where $\bar{\mathfrak{H}}$ is conjugate 
linearly isomorphic to $\mathfrak{H}$ (the image of $\xi$ is written 
$\bar{\xi}$), and the actions are defined by $n \bar{\xi} m = \overline{m^* 
\xi n^*}.$
\end{definition}

\begin{lemma} \label{T:inverse}
Suppose $L(\mathcal{M})$ and $R(\mathcal{N})$ are commutants on $\mathfrak{H}$.  Then
$${}_\mathcal{N} \bar{\mathfrak{H}}_\mathcal{M} \otimes_\mathcal{M} 
{}_\mathcal{M} \mathfrak{H}_\mathcal{N} \simeq {}_\mathcal{N} 
L^2(\mathcal{N})_\mathcal{N}.$$
\end{lemma}

\begin{proof}
If ${}_\mathcal{M} \mathfrak{H} \simeq {}_\mathcal{M} R^2(\mathcal{M})q$, 
then $\mathcal{N} \simeq qM_\infty(\mathcal{M})q$ by \eqref{E:rcomm}, and 
$\bar{\mathfrak{H}}_\mathcal{M} \simeq qC^2(\mathcal{M})_\mathcal{M}.$  By 
Proposition \ref{T:colxrow} and the comment preceding Proposition \ref{T:module},
$${}_\mathcal{N} \bar{\mathfrak{H}}_\mathcal{M} \otimes_\mathcal{M} 
{}_\mathcal{M} \mathfrak{H}_\mathcal{N} \simeq {}_\mathcal{N} 
(qL^2(M_\infty(\mathcal{M}))q)_\mathcal{N} \simeq {}_\mathcal{N} L^2(q 
M_\infty(\mathcal{M})q)_\mathcal{N} \simeq {}_\mathcal{N} L^2(\mathcal{N}) 
_\mathcal{N}.$$
\end{proof}

Lemma \ref{T:inverse} was first proven by Sauvageot (in another way).  In our situation it means

$$F_{\bar{\mathfrak{H}}} \circ F_\mathfrak{H} ({}_\mathcal{N} \mathfrak{K}) \simeq L^2(\mathcal{N}) \otimes_\mathcal{N} {}_\mathcal{N} \mathfrak{K} \simeq {}_\mathcal{N} \mathfrak{K}.$$
(Here we have used the associativity of the RTP, which is most easily seen from the explicit construction in Section \ref{S:analytic}.)  We conclude that $F_{\bar{\mathfrak{H}}} \circ F_\mathfrak{H}$ is equivalent to the identity functor on $Left \:L^2 Mod(\mathcal{N})$, and by a symmetric argument $F_\mathfrak{H} \circ F_{\bar{\mathfrak{H}}}$ is equivalent to the identity functor on $Left \:L^2 Mod(\mathcal{M})$.  Thus $F_\mathfrak{H}$ is a Morita equivalence, and the first implication of Theorem \ref{T:morita} is proved.

Now let $F$ be a Morita equivalence from $Left\:L^2 Mod (\mathcal{N})$ to $Left\:L^2 Mod(\mathcal{M})$.  Then $F$ must take $R^2(\mathcal{N})$ to a module isomorphic to $R^2(\mathcal{M})$, because each is in the unique isomorphism class which absorbs all other modules.  (This is meant in the sense that ${}_\mathcal{N} R^2(\mathcal{N}) \oplus {}_\mathcal{N} \mathfrak{H} \simeq {}_\mathcal{N} R^2(\mathcal{N})$; $R^2(\mathcal{N})$ is the ``infinite element" in the monoid $V(M_\infty(\mathcal{N}))$.)  Being an invertible *-functor, $F$ implements a *-isomorphism of commutants - call it $\sigma$, not $F$, to ease the notation:

\begin{equation} \label{E:stab}
\sigma: M_\infty(\mathcal{N}) \overset{\sim}{\to} M_\infty (\mathcal{M}).
\end{equation}
Apparently we have

\begin{equation} \label{E:morita}
F(R^2(\mathcal{N}) q) \simeq R^2(\mathcal{M}) \sigma(q). 
\end{equation}

Before continuing the argument, we need an observation: isomorphic algebras have isomorphic standard forms.  Specifically, we may write $L^2(M_\infty(\mathcal{N}))$ as the GNS construction for $\varphi \in M_\infty(\mathcal{N})_*^+$ and obtain the isomorphism
$$(\sigma^{-1})^\text{t}: L^2(M_\infty(N)) \overset{\sim}{\to} L^2(M_\infty(M)), \qquad (\sigma^{-1})^\text{t}: x \varphi^{1/2} \mapsto \sigma(x) (\varphi \circ \sigma^{-1})^{1/2}.$$ 
Note that $(\sigma^{-1})^\text{t}(x \xi y) = \sigma(x) [(\sigma^{-1})^\text{t}(\xi)] \sigma(y).$

Now consider the RTP functor for the bimodule 

$${}_\mathcal{M} \mathfrak{H}_\mathcal{N} = {}_{\sigma^{-1}(\mathcal{M})} \sigma^{-1}(e_{11}^\mathcal{M}) C^2(\mathcal{N})_\mathcal{N}.$$
Its action is

$$R^2(\mathcal{N})q \mapsto {}_{\sigma^{-1}(\mathcal{M})} \sigma^{-1}(e_{11}^\mathcal{M}) L^2(M_\infty(\mathcal{N})) q \overset{(\sigma^{-1})^\text{t}}{\simeq} {}_\mathcal{M} e_{11}^\mathcal{M} L^2(M_\infty(\mathcal{M}))\sigma(q)$$
$$\simeq {}_\mathcal{M} R^2(\mathcal{M}) \sigma(q) \simeq F(R^2(\mathcal{N})q).$$
We conclude that $F$ is equivalent to $F_\mathfrak{H}$, which finishes the proof of Theorem \ref{T:morita}.

Notice that \eqref{E:stab} and \eqref{E:morita} can also be used to define a functor; this gives us

\begin{corollary}
For two von Neumann algebras $\mathcal{M}$ and $\mathcal{N}$, the following are equivalent:

\begin{enumerate}

\item $\mathcal{M}$ and $\mathcal{N}$ are Morita equivalent;

\item $M_\infty(\mathcal{N}) \simeq M_\infty (\mathcal{M})$;

\item there is a bimodule ${}_\mathcal{M} \mathfrak{H}_\mathcal{N}$ where the actions are commutants of each other;

\item there is a projection $q \in M_\infty(\mathcal{M})$ with central support 1 such that
$$q M_\infty(\mathcal{M})q \simeq \mathcal{N}.$$  
\end{enumerate}

\end{corollary}
(The \textit{central support} of $x \in \mathcal{M}$ is the least projection $z$ in the center of $\mathcal{M}$ satisfying $x=zx$.)

\textbf{Example continued.} $M_3$ and $M_5$ are Morita equivalent.  This can be deduced easily from condition (2) or (4) of the corollary above.  It also follows from the (Hilbert) equivalence bimodule ${}_{M_3} HS(M_{3\times 5})_{M_5}$, where ``$HS$" indicates the Hilbert-Schmidt norm; this bimodule gives us an RTP functor which is a Morita equivalence.  Regardless of the construction, the equivalence will send (functorially) $n$ copies of $\mathbb{C}^5$ to $n$ copies of $\mathbb{C}^3$.  Apparently Morita equivalence is a coarse relation on type I algebras - it only classifies the \textit{center} of the algebra (up to isomorphism).  At the other extreme, Morita equivalence for type III algebras is the same as algebraic isomorphism; Morita equivalence for type II algebras is somewhere in the middle ([R], Sec. 8).

\bigskip

For a bimodule ${}_\mathcal{M} \mathfrak{H}_\mathcal{N}$ where the algebras are not necessarily commutants, the functor \eqref{E:rtpfunctor} still makes sense.  To get a more tractable object, we may consider the domain and range to be isomorphism classes of modules:

\begin{equation} \label{E:bmorph}
\pi_\mathfrak{H}: V(M_\infty(\mathcal{N})) \to V(M_\infty(\mathcal{M}));
\end{equation}
$$F_\mathfrak{H}(R^2(\mathcal{N})q) = {}_\mathcal{M} \mathfrak{H}_\mathcal{N} \otimes_\mathcal{N} R^2(\mathcal{N})q \simeq R^2(\mathcal{M}) \pi_\mathfrak{H}([q]).$$

We call this the \textit{bimodule morphism} associated to $\mathfrak{H}$, a sort of ``skeleton" for the correspondence.  It follows from Proposition \ref{T:corresp} that if the bimodule is $\rho: \mathcal{N} \hookrightarrow qM_\infty(\mathcal{M})q$, then $\pi_\mathfrak{H}$ is nothing but $\rho^\infty$, the amplification of $\rho$ to $M_\infty(\mathcal{N})$, restricted to equivalence classes of projections.

This has an important application to inclusions of algebras.  We have seen in Section \ref{S:algebraic} that a unital inclusion $\mathcal{N} \overset{\rho}{\subset} \mathcal{M}$ is equivalent to a bimodule ${}_\mathcal{M} L^2(\mathcal{M})_\mathcal{N}$.  When the correspondence $\rho$ is surjective, the module generates a Morita equivalence via its RTP functor, and the induced bimodule morphism is an isomorphism of monoids.  When $\mathcal{N} \ne \mathcal{M}$, it is natural to expect that the bimodule morphism gives us a way to measure the relative size, or \textit{index}, of $\mathcal{N}$ in $\mathcal{M}$.  (For readers unfamiliar with this concept, the index of an inclusion $\mathcal{N} \subset \mathcal{M}$ is denoted $[\mathcal{M}: \mathcal{N}]$ and is analogous to the index of a subgroup.  The kernel of this idea goes back to Murray and von Neumann, but the startling results of Jones [J] in the early 1980's touched off a new wave of investigation.  We recommend the exposition [K] as a nice starting point.)

For algebras of type I or II, the index can be calculated in terms of bimodule morphisms.  (There are also broader definitions of index which require a conditional expectation (=norm-decreasing projection) from $\mathcal{M}$ onto $\mathcal{N}$.)  This amounts largely to rephrasing and extension of the paper [Jol], and we do not give details here.  Very briefly, let $\pi: V(M_\infty(\mathcal{M})) \to V(M_\infty(\mathcal{M}))$ be the bimodule morphism for

\begin{equation} \label{E:sym}
({}_\mathcal{M} L^2(\mathcal{M})_\mathcal{N}) \otimes_\mathcal{N} ({}_\mathcal{N} L^2(\mathcal{M})_\mathcal{M}).
\end{equation}
When $\mathcal{M}$ is a factor, $\pi$ is a monoidal morphism on the extended nonnegative integers (type I) or extended nonnegative reals (type II).  It must be multiplication by a scalar, and this scalar is the index.  If $\mathcal{M}$ is not a factor, the index is the spectral radius of $\pi$, provided that $V(M_\infty(\mathcal{M}))$ is endowed with some extra structure (at present it is not even a vector space).

\textbf{Example.} Consider the correspondence

$${}_{M_6} L^2(M_6)_{M_3}; \qquad M_3 \simeq M_3 \otimes I \subset M_3 \otimes M_2 \simeq M_6.$$
The image of ${}_{M_6}L^2(M_6)$ under the RTP functor for \eqref{E:sym} is

$$\left({}_{M_6} L^2(M_6)_{M_3} \right) \otimes_{M_3} \left({}_{M_3} L^2(M_6)_{M_6} \right) \otimes_{M_6} \left({}_{M_6}L^2(M_6)\right)$$

$$\simeq {}_{M_6} L^2(M_6)_{M_3} \otimes_{M_3} {}_{M_3} L^2(M_6)$$
(now counting the dimensions of the Hilbert spaces)

$$\simeq {}_{M_6} HS(M_{12\times3})_{M_3} \otimes_{M_3} {}_{M_3} HS(M_{3\times12}) \simeq {}_{M_6} HS(M_{12\times12}) \simeq {}_{M_6} HS(M_{6\times24}).$$
We have gone from 6 copies of $\mathbb{C}^6$ to 24 copies; that is,

$$\pi: (\mathbb{Z}_+ \cup \infty) \to (\mathbb{Z}_+ \cup \infty), \qquad 6 \mapsto 24.$$
Apparently the index is 4, which is also the ratio of the dimensions of the algebras.  

\section{Analytic approaches to RTPs} \label{S:analytic}

As mentioned in Section \ref{S:algebraic}, we cannot expect the expression $\xi \otimes_\mathcal{M} \eta$ to make sense for every pair of vectors $\xi, \eta$.  In essence, the problem is that the product of two $L^2$ vectors is $L^1$, and an $L^1$ space does not lie inside its corresponding $L^2$ space unless the underlying measure is atomic.  Densities add, even in the noncommutative setting, and so the product in \eqref{E:l2l2} ``should" be an $L^1$ matrix.  To make this work at the vector level, we need to decrease the density by $1/2$ without affecting the ``outside" action of the commutants... and the solution by Connes and Sauvageot ([Sa], [C2]) is almost obvious: choose a faithful $\varphi \in \mathcal{M}_*^+$ and put $\varphi^{-1/2}$ in the middle of the product.  That is,

\begin{equation} \label{E:neghalf}
\xi \otimes_\varphi \eta = (\xi \varphi^{-1/2}) \eta.
\end{equation}
This requires some explanation.

Fix faithful $\varphi \in \mathcal{M}_*^+$ and row and column representations of $\mathfrak{H}$ and $\mathfrak{K}$ as in \eqref{E:rowrep}.  We define
$$\mathcal{D}(\mathfrak{H}, \varphi) = \left\{ \left(
\begin{smallmatrix}
x_1 \varphi^{1/2} \\
x_2 \varphi^{1/2} \\
\vdots
\end{smallmatrix}
\right) \in \mathfrak{H} \: : \: \sum x_n^* x_n \text{ exists in }\mathcal{M} \right\}.$$
$\mathcal{D}(\mathfrak{H}, \varphi)$ is dense in $\mathfrak{H}$, and its elements are called \textit{$\varphi$-left bounded} vectors [C1]. 

Now by \eqref{E:neghalf} we mean the following: for $\xi \in \mathcal{D}(\mathfrak{H}, \varphi)$, we simply erase the symbol $\varphi^{1/2}$ from the right of each entry, then carry out the multiplication.  The natural domain is $\mathcal{D}(\mathfrak{H}, \varphi) \times \mathfrak{K}$.  Visually,

\begin{equation} \label{E:visual}
\left(
\left(
\begin{smallmatrix}
x_1 \varphi^{1/2} \\
x_2 \varphi^{1/2} \\
\vdots
\end{smallmatrix}
\right),
\left(
\begin{smallmatrix}
\eta_1 & \eta_2 & \cdots
\end{smallmatrix}
\right) \right) \mapsto \left(
\begin{smallmatrix}
x_1 \varphi^{1/2} \\
x_2 \varphi^{1/2} \\
\vdots
\end{smallmatrix}
\right)
(\begin{smallmatrix} \varphi^{-1/2} \end{smallmatrix})
\left(
\begin{smallmatrix}
\eta_1 & \eta_2 & \cdots
\end{smallmatrix}
\right) = (x_i \eta_j).
\end{equation}
For $\varphi \ne \psi \in \mathcal{M}_*^+,$ we cannot expect $\xi \otimes_\varphi \eta = \xi \otimes_\psi \eta$ even if both are defined, although the reader familiar with modular theory will see that 

\begin{equation} \label{E:convert}
\xi \otimes_\varphi \eta = (\xi \varphi^{-1/2}) \eta = (\xi \varphi^{-1/2} \psi^{1/2} \psi^{-1/2}) \eta = \xi (D\varphi: D\psi)_{i/2} \otimes_\psi \eta.
\end{equation}
(An interpretation of the symbols $\varphi^{1/2}, \varphi^{-1/2}$ as unbounded operators will be discussed in the next section.)

Now we define $\mathfrak{H} \otimes_\varphi \mathfrak{K}$ to be the closed linear span of the vectors $\xi \otimes_\varphi \eta$ inside $L^2(M_\infty(\mathcal{M}))$.  Up to isomorphism, this is independent of $\varphi$.  (We know this because of functoriality; the ``change of weight" isomorphism is densely defined by \eqref{E:convert}.)

The given definition for $\mathcal{D}(\mathfrak{H}, \varphi) \subset \mathfrak{H}$ makes it seem dependent on the choice of column representation.  That this is not so can be seen by noting (as in \eqref{E:rightiso}) that the intertwining isomorphism is given by $L(v)$ for some partial isometry $v \in M_\infty(\mathcal{M}).$  But let us also mention a method of defining the same RTP construction without representing $\mathfrak{H}$ and $\mathfrak{K}$.  First notice that $\mathcal{D}(\mathfrak{H}, \varphi)$ can also be defined as the set of vectors $\xi$ for which
$$\pi_\ell^\varphi(\xi): L^2(\mathcal{M})_\mathcal{M} \to \mathfrak{H}_\mathcal{M}, \qquad \varphi^{1/2} x \mapsto \xi x,$$
is bounded.  (A more suggestive (and rigorous) notation would be $L(\xi \varphi^{-1/2})$.)  Now we consider an inner product on the algebraic tensor product $\mathcal{D}(\mathfrak{H}, \varphi) \otimes \mathfrak{K}$, defined on simple tensors by

\begin{equation} \label{E:innprod}
<\xi_1 \otimes_\varphi \eta_1 \mid \xi_2 \otimes_\varphi \eta_2> = <\pi_\ell^\varphi(\xi_2)^* \pi_\ell^\varphi(\xi_1) \eta_1 \mid \eta_2>.
\end{equation}
The important point here is that $\pi_\ell^\varphi(\xi_2)^* \pi_\ell^\varphi(\xi_1)$ $\in \mathcal{L}(L^2(\mathcal{M})_\mathcal{M}) = \mathcal{M}$.  The closure of $\mathcal{D}(\mathfrak{H}, \varphi) \otimes \mathfrak{K}$ in this inner product, modulo the null space, is once again $\mathfrak{H} \otimes_\varphi \mathfrak{K}$.

(If we \textit{do} choose a row representation as in \eqref{E:visual}, we have

$$ \pi_\ell^\varphi \left( \left(
\begin{smallmatrix}
x_1 \varphi^{1/2} \\
x_2 \varphi^{1/2} \\
\vdots
\end{smallmatrix}
\right) \right) = L \left( \left(
\begin{smallmatrix}
x_1 \\
x_2 \\
\vdots
\end{smallmatrix}
\right) \right).)$$

The paper [F] contains more exposition of this approach, including some alternate constructions.

\section{Realization of the relative tensor product as composition of unbounded operators} \label{S:modalg}

In this section we briefly indicate how \eqref{E:neghalf} can be rigorously justified and extended.  Readers are referred to the sources for all details.

In his pioneering theory of noncommutative $L^p$ spaces, Haagerup [H] established a linear isomorphism between $\mathcal{M}_*^+$ and a class of positive unbounded operators affiliated with the core of $\mathcal{M}$.  (The core, well-defined up to isomorphism, is the crossed product of $\mathcal{M}$ with one of its modular automorphism groups.)  We will denote the operator corresponding to the positive functional $\varphi$ by $\varphi$ also.  These operators are $\tau$-measurable (see the next section), where $\tau$ is the canonical trace on the core, and so they generate a certain graded *-algebra: positive elements of $L^p(\mathcal{M})$ are defined to be operators of the form $\varphi^{1/p}$.  The basic development of this theory can be found in [Te]; our choice of notation is influenced by [Y], where it is called a \textit{modular algebra}.

The composition of two $L^2$ operators is an $L^1$ operator, and it turns out that \eqref{E:neghalf} can be rigorously justified [S] as an operator equation.  (This is not automatic, as operators like $\varphi^{-1/2}$ are not $\tau$-measurable and require more delicate arguments.)  In fact, there is nothing sacred about half-densities.  With the recent development of noncommutative $L^p$ modules [JS], one can allow relative tensor products to be bifunctors on $Right \: L^p(\mathcal{M}) \times Left \: L^q (\mathcal{M})$, with range in a certain $L^r$ space.  The mapping is densely-defined by

$$\xi \otimes_\varphi \eta \triangleq (\xi \varphi^{\frac{1}{r} - \frac{1}{p} - \frac{1}{q}}) \eta.$$

In the case of an RTP of $L^\infty$ modules (or more generally, Hilbert C*-modules), the middle term is trivial and there is no change of density.  This explains why there are no domain issues in defining a vector-valued RTP of Hilbert C*-modules [R].

Let us mention that the recent theory of operator bimodules, in which vectors can be realized as \textit{bounded} operators, allows a variety of relative tensor products over C*-algebras [AP].  This can be naturally viewed as a generalization of the theory of Banach space tensor products, which corresponds to a C*-algebra of scalars.

\section{Preclosedness of the map $(\xi, \eta) \mapsto \xi \otimes_\varphi \eta$} \label{S:preclosed}

Our purpose in this final section is to study when the relative tensor map is \textit{preclosed}.  This is a weaker condition than that of Falcone, who studied (effectively) when the map was \textit{bounded}.  We begin with a base case: a factor, two standard modules, and a simple product.  With the usual notation $\mathcal{B}_\varphi$ for $\mathcal{D}(L^2(\mathcal{M}), \varphi)$, the relevant map is 

$$\mathcal{B}_\varphi \times L^2(\mathcal{M}) \ni (\xi, \eta) \mapsto \xi \otimes_\varphi \eta \in L^2(\mathcal{M}).$$
This is bilinear: we take ``preclosed" to mean that if $\xi_n \to \xi \in \mathcal{B}_\varphi,$ $\eta_n \to \eta$, $\xi_n \otimes_\varphi \eta_n \to \zeta$, then necessarily $\zeta = \xi \otimes_\varphi \eta.$  We will also consider several variations: changing the domain to an algebraic tensor product, allowing non-factors, and allowing arbitrary modules.

Readers unfamiliar with von Neumann algebras will find this section more technical, and any background we can offer here is sure to be insufficient.  Still, we introduce the necessary concepts in hopes that the non-expert will at least find the statements of the theorems accessible.

\bigskip

A \textit{weight} is an ``unbounded positive linear functional": a linear map from $\mathcal{M}_+$ to $[0, +\infty]$.  We will always assume that weights are \textit{normal}, so $x_\alpha \nearrow x$ strongly $\Rightarrow \varphi(x_\alpha) \nearrow \varphi(x)$; and \textit{semifinite}, so $\{x \in \mathcal{M}_+ \mid \varphi(x) < \infty \}$ is $\sigma$-weakly dense in $\mathcal{M}_+$.  We can still define RTPs for faithful weights, but now $\mathcal{B}_\varphi = \{x\varphi^{1/2} \mid \varphi(x^* x) < \infty\} \subset L^2(\mathcal{M})$.  For details of the representations associated to weights, see [T2].

A weight $\tau$ which satisfies $\tau(xy)= \tau(yx)$ on its domain of definition will be called a \textit{trace} (more properly called a ``tracial weight").  An algebra which admits a faithful trace $\tau$ is \textit{semifinite}; if in addition we can have $\tau(1) < \infty$, it is \textit{finite}.  This facilitates the following classification of factors: a factor with $n$ orthogonal minimal projections is type $\text{I}_n$ (possibly $n=\infty$); a semifinite factor without minimal projections is type $\text{II}_1$ if finite and type $\text{II}_\infty$ if not; a factor which is not semifinite is type III.  The reader will note that this refines our previous definitions of type, as a trace is exactly the object which orders the equivalence classes of projections.  Obviously, there is much more to be said, and most of it can be found in [T1].

For a faithful trace $\tau$ on semifinite $\mathcal{M}$, it is useful to consider the \textit{$\tau$-measure topology} [N].  This is a uniform topology with neighborhoods of 0 given by 

$$N(\delta, \epsilon) = \{x \in \mathcal{M} \mid \exists p \in \mathcal{P}(\mathcal{M}) \text{ with } \tau(p^\perp) < \delta, \|xp\| < \epsilon\}.$$
The closure of $\mathcal{M}$ in this topology can be identified as a space of closed, densely-defined operators affiliated with $\mathcal{M}$.  It is denoted $\mathfrak{M}(\mathcal{M})$ and actually forms a *-algebra to which $\tau$ extends naturally.  (The $\tau$-measurability of an operator $T$ is equivalent to the assertion that $\tau(e(\lambda)^\perp) < \infty$ for some spectral projection $e(\lambda)$ of $|T|$, so we get that $\mathfrak{M}(\mathcal{M}) = \mathcal{M}$ if $\mathcal{M}$ is atomic.)  It follows from modular theory that every weight on $(\mathcal{M}, \tau)$ is of the form $\tau_h = $``$\tau(h \cdot)$" for some closed, densely-defined, and positive operator $h$.  In case $h$ is not $\tau$-measurable, this is to be interpreted as

$$\lim_{\varepsilon \to 0} \tau(h_\varepsilon^{1/2} \cdot h_\varepsilon^{1/2}), \text{ where } h_\varepsilon = h(1 + \varepsilon h)^{-1}.$$

Finally, the presence of a faithful trace $\tau$ allows us to introduce the spaces

$$L^p(\mathcal{M}, \tau) = \{T \in \mathfrak{M}(\mathcal{M}) \mid \tau(|T|^p) = \|T\|^p < \infty \},$$
which are antecedent to Haagerup's.  Exposition can be found in [N].  Here we will only need $L^2(\mathcal{M}, \tau)$, which is a standard form and in particular isomorphic as a left module to any faithful GNS representation $\mathfrak{H}_\varphi$.  It is easy to check that the norm topology in $L^2(\mathcal{M}, \tau)$ is stronger than the $\tau$-measure topology.

\begin{theorem} \label{T:preclosed}
Let $\mathcal{M}$ be a factor.  The map
\begin{equation}
\mathfrak{B}_\varphi \times L^2(\mathcal{M}) \to L^2(\mathcal{M}): \quad (\xi,\eta) \mapsto \xi \otimes_\varphi \eta
\end{equation}
is preclosed iff $\mathcal{M} = (\mathcal{M},\tau)$ is semifinite and $h^{-1}$ is $\tau$-measurable, where $\varphi = \tau_h$.
\end{theorem}

\begin{proof}
The proof is by consideration of cases.  

\bigskip

\textit{$\mathcal{M}$ is type III:}  Choose a projection $e_0$ so $\varphi(e_0)=c<\infty$.  Find orthogonal projections $f_1$, $g_1$ with $e_0=f_1+g_1$.  Set $e_1 = f_1$ if $\varphi(f_1) \le \varphi(g_1)$ and $e_1=g_1$ otherwise.  Continuing in this fashion gives a sequence of projections, necessarily $\sigma$-finite, with $\varphi(e_n) \le c/{2^n}$.  Thus all the $e_n$ are Murray-von Neumann equivalent, and there exist partial isometries $v_n$ with $v_n^* v_n = e_0$, $v_n v_n^* = e_n$.  Implementing multiplication,

$$(nv_n^*\varphi^{1/2}, (1/n) v_n \varphi^{1/2}) \mapsto v_n^*v_n \varphi^{1/2} = e_0 \varphi^{1/2}.$$
But $$||nv_n^*\varphi^{1/2}||^2 = n^2(\varphi(e_n)) \to 0;$$

$$||(1/n) v_n \varphi^{1/2}||^2 = (1/n^2) (\varphi(e_0)) \to 0;$$
thus the multiplication is not preclosed.

\bigskip

\textit{$\mathcal{M} = (\mathcal{M},\tau)$ is semifinite and $h^{-1}$ is not $\tau$-measurable:}  First note that the measurability of $h^{-1}$ does not depend on the choice of $\tau$.  Writing $h = \int \lambda de(\lambda)$, the hypothesis is that $\tau(e(\lambda)) = \infty,\;\forall \lambda$.  Choose a projection $e_0$ with $\varphi(e_0) < \infty$ and $\tau(e_0) < \infty$.  Then $e(1/{n^3})$ has a subprojection $e_n$ which is equivalent to $e_0$.  The above construction again shows that the map is not preclosed, except that

$$||nv_n^*\varphi^{1/2}||^2 = n^2(\varphi(e_n)) = n^2(\tau(he_n)) \le (1/n) \tau(e_n) \to 0.$$

\bigskip

\textit{$\mathcal{M} = (\mathcal{M},\tau)$ is semifinite and $h^{-1}$ is $\tau$-measurable:}  We assume

\begin{equation} \label{E:conv}
x_n \varphi^{1/2} \to x \varphi^{1/2}, \quad \eta_n \to \eta, \quad x_n \eta_n \to \zeta,
\end{equation}
and want to show $\zeta = x \eta$.  Set

$$\mathfrak{n} = \{x \in \mathcal{M} \mid \overline{x h^{1/2}} \in L^2(\mathcal{M}, \tau) \}; \qquad \mathfrak{n}_\varphi = \{x \in \mathcal{M} \mid \varphi(x^* x) < \infty \},$$
both of which are strongly dense in $\mathcal{M}.$  (The bar stands for graph closure.)  Then

$$\pi: \mathfrak{n} \varphi^{1/2} \to L^2(\mathcal{M}, \tau); \quad x\varphi^{1/2} \mapsto \overline{x h^{1/2}}$$
densely defines a left module Hilbert space isomorphism from $\mathfrak{H}_\varphi$ to $L^2(\mathcal{M}, \tau)$; denote its extension by $\pi$ as well.  Recalling that $h^{-1/2}$ is $\tau$-measurable by assumption,

$$\rho: \mathfrak{n}_\varphi \to \mathfrak{M}(\mathcal{M}); \quad x \mapsto \pi(x\varphi^{1/2}) h^{-1/2}$$
is well-defined and the identity map on $\mathfrak{n}$.  It is also strong-measure continuous:

$$x_\alpha \overset{s}{\to} x \Rightarrow x_\alpha \varphi^{1/2} \overset{\mathfrak{H}_\varphi}{\to} x \varphi^{1/2} \Rightarrow \pi(x_\alpha \varphi^{1/2}) \overset{L^2}{\to} \pi(x \varphi^{1/2})$$ $$\Rightarrow \pi(x_\alpha \varphi^{1/2}) \overset{m}{\to} \pi(x \varphi^{1/2}) \Rightarrow \pi(x_\alpha \varphi^{1/2})h^{-1/2} \overset{m}{\to} \pi(x \varphi^{1/2})h^{-1/2},$$
where we used that multiplication is jointly continuous in the measure topology.  We may conclude that $\rho$ is the identity on all of $\mathfrak{n}_\varphi$.

Implementing the isomorphism $\pi$, \eqref{E:conv} becomes

\begin{equation} \label{E:convergences}
\pi(x_n \varphi^{1/2}) \to \pi(x \varphi^{1/2}), \quad \pi(\eta_n) \to \pi(\eta), \quad x_n \pi(\eta_n) \to \pi(\zeta).
\end{equation}
The convergences in \eqref{E:convergences} are also in measure; by the foregoing discussion we have

$$x_n \pi(\eta_n) = \pi(x_n \varphi^{1/2}) h^{-1/2} \pi(\eta_n) \to \pi(x \varphi^{1/2}) h^{-1/2} \pi(\eta) = x \pi(\eta)$$
in measure as well.

The measure topology is also Hausdorff, so $\pi(\zeta) = x \pi(\eta)$ and therefore $\zeta = x \eta$.

\end{proof}

The map $\rho$ suggests a schematic recovery of the ``operators" in $\mathfrak{H}_\varphi$:

\begin{equation} \label{E:schematic}
\{\pi_\ell^\varphi(\xi)\mid \xi \in \mathfrak{H}_\varphi \} = L^2(\mathcal{M}, \tau) h^{-1/2}.
\end{equation}
Such operators are densely-defined but in general not closable (or may have multiple closed extensions [Sk]).  Not surprisingly, then, the right-hand side of \eqref{E:schematic} may be only formal.  The condition on $h$ in Theorem \ref{T:preclosed} makes the equality \eqref{E:schematic} rigorous, as the products on the right-hand side are well-defined $\tau$-measurable operators.  Note that $h$ and $h^{-1}$ are automatically $\tau$-measurable when $\mathcal{M}$ is finite, and in this case all multiplications and isomorphisms between GNS representations stay within $\mathfrak{M}(\mathcal{M})$, and all operators are closed - a version, somewhat oblique, of the $T$-theorem of Murray and von Neumann.

\begin{theorem} \label{T:precten}

Let $\mathcal{M}$ be a factor, and consider $\mathfrak{B}_\varphi \otimes_{\text{alg}} L^2(\mathcal{M})$ as a subspace of the Hilbert space tensor product $L^2(\mathcal{M}) \otimes L^2(\mathcal{M})$.  The linear map

\begin{equation}
\mathfrak{B}_\varphi \otimes_{\text{alg}} L^2(\mathcal{M}) \to L^2(\mathcal{M}): \quad \sum \xi_n \otimes \eta_n \mapsto \sum \xi_n \otimes_\varphi \eta_n
\end{equation}
is preclosed iff $\mathcal{M} = (\mathcal{M},\tau)$ is atomic and $\tau(h^{-1})<\infty$, where $\varphi = \tau_h$.  In this case it is actually a bounded map, with norm $\tau(h^{-1})^{1/2}$.

\end{theorem}

\begin{proof}
When $\mathcal{M}$ is type III, the map is not preclosed by the previous theorem.  We will therefore fix a trace $\tau$, set $\varphi = \tau_h$, use the decomposition $h=\int \lambda \,de(\lambda)$, and view all vectors as elements of $L^2(\mathcal{M}, \tau)$.  (When $\mathcal{M}$ is type I, we assume that $\tau$ is normalized so that $\tau(p)=1$ for any minimal projection $p$.)  The rest of the proof is again by cases.

\bigskip

\textit{$\mathcal{M}$ is type II:}  Choose $p<e(\lambda)$ for some $\lambda$ with $\tau(p)=c<\infty$.  For each $k$, break up $p$ into equivalent orthogonal projections as $\sum_{n=1}^k p_n^k.$  Consider the tensors

$$T_k = \sum p_n^k h^{1/2} \otimes p_n^k \mapsto \sum p_n^k = p.$$
Since the $p_n^k$ are orthogonal,

$$||T_k||^2 = \sum \tau(p_n^k h) \tau(p_n^k) \le \sum \left(\frac{\lambda c}{k}\right)\left(\frac{c}{k}\right) = \frac{\lambda c^2}{k} \to 0$$
and the map is not preclosed.

\bigskip

\textit{$\mathcal{M}$ is type $I_\infty$ and $\tau(e(\lambda)) =\infty$ for some $\lambda$:}  Fix an orthogonal sequence of minimal projections $\{p_n\}$, $p_n < e(\lambda)$.  The equivalence gives partial isometries with $v_n^*v_n = p_1$, $v_n v_n^* = p_n$.  Then

$$T_k = \sum_{n=1}^k \frac{1}{k}(v_n^*h^{1/2} \otimes v_n) \mapsto \sum \frac{1}{k} v_n^*v_n = p_1;$$

$$ ||T_k||^2 = \frac{1}{k^2} \sum \tau (p_n h) \tau (p_1) \le \frac{1}{k^2} \sum \lambda = \frac{\lambda}{k} \to 0$$
and the map is not preclosed.

\bigskip

In the only remaining situation, $\mathcal{M}$ is type $\text{I}$ and $h$ is diagonalizable.  Let $\{\lambda_n\}$ be the eigenvalues (with repetition), arranged in nondecreasing order.  We will write all matrices with respect to the basis of eigenvectors.

\bigskip

\textit{If $s_k = \sum_{n=1}^k \frac{1}{\lambda_n} \nearrow \infty$:}  Consider

\begin{equation} \label{E:nonnuclear}
T_k = \sum_{n=1}^k \frac{1}{s_k \lambda_n} e_{1n} h^{1/2} \otimes e_{n1} \mapsto \frac{1}{s_k} \sum \frac{1}{\lambda_n} e_{11} = e_{11};
\end{equation}

$$||T_k||^2 = \sum \frac{1}{s_k^2 \lambda_n^2} \tau(e_{nn}h) \tau(e_{11}) = \frac{1}{s_k^2} \sum \frac{1}{\lambda_n} = \frac{1}{s_k} \to 0$$
and the map is not preclosed.

\bigskip

\textit{If $s_k = \sum_{n=1}^k \frac{1}{\lambda_n} \nearrow C < \infty$; that is, $\tau(h^{-1}) < \infty$:}  We show that the map is bounded on finite tensors of the form $T = \sum_{ij} e_{ij} \otimes y^{ij}.$  We have

$$T \mapsto S = \sum_{ij} e_{ij} h^{-1/2} y^{ij} = \sum_{ij} e_{ij} \lambda_j^{-1/2} y^{ij} = \sum_{ik} \left( \sum_j \lambda_j^{-1/2} y^{ij}_{jk} e_{ik} \right).$$

By Cauchy-Schwarz,

$$\|S\|^2 = \sum_{ik} \left| \sum_j \lambda_j^{-1/2} y^{ij}_{jk} \right|^2 \le \sum_{ik} \left( \left[\sum_j \lambda_j^{-1} \right] \left[ \sum_j |y^{ij}_{jk}|^2 \right]\right)$$ $$\le C \sum_{ijk} | y^{ij}_{jk}|^2 \le C \sum_{ijkl} | y^{ij}_{lk}|^2 = C \|T\|^2.$$

Since such tensors are dense in the Hilbert space tensor product, we may conclude that the norm of the map is $\le C^{1/2}.$  But the tensors $T_k$ from \eqref{E:nonnuclear} show that the norm is at least $C^{1/2}$.
\end{proof}

We now extend Theorem \ref{T:preclosed} to the non-factor case.  A general von Neumann algebra is a direct integral of factors (see [T1] for details), and weights on the algebra disintegrate as well.

\begin{proposition} \label{T:genprec}
Let $\mathcal{M}$ be a von Neumann algebra with central decomposition $\int^\oplus_\Gamma \mathcal{M}(\omega) d\mu(\omega).$  The map
\begin{equation}
\mathfrak{B}_\varphi \times L^2(\mathcal{M}) \to L^2(\mathcal{M}): \quad (\xi,\eta) \mapsto \xi \otimes_\varphi \eta
\end{equation}
is preclosed iff $\mathcal{M} = (\mathcal{M},\tau)$ is semifinite and \\$(\star)$ $h(\omega)^{-1}$ is $\tau(\omega)$-measurable for $\mu$-a.e. $\omega$, where $\varphi = \tau_h$.

\end{proposition}

\begin{proof}
If $\mathcal{M}$ contains a summand of type III, the construction from Theorem \ref{T:preclosed}, with the added restriction that $f_n$ and $g_n$ are chosen with equal central support, demonstrates that the map is not preclosed.

If there is a trace $\tau$ for which $\varphi = \tau_h$ and $h^{-1}$ is $\tau$-measurable, then the argument in Theorem \ref{T:preclosed} still shows that the map is preclosed.  We will see that this possibility is equivalent to $(\star)$.

First note that $(\star)$ is independent of the trace chosen, as the choice of a different trace changes a.e. $h(\omega)$ by a constant factor.  If $(\star)$ does not hold, fix any trace $\tau$, write $\varphi = \tau_h$, and let $\{e(\lambda)\}$ be the spectral projections of $h$.  By hypothesis, we can find a nonzero central projection $z$ with $z e(\lambda)$ a properly infinite projection for all $\lambda$.  The second construction of Theorem \ref{T:preclosed} shows that the map is not preclosed, where we choose all $e_n$, including $e_0$, with central support $z$.

Now suppose that $(\star)$ holds.  We may choose a trace $\tau$ which factors as $\tilde{\tau} \circ \Phi$, where $\Phi$ is an extended center-valued trace and $\tilde{\tau}$ is a trace on the center with $\tilde{\tau}(1) < \infty.$  Let $h$ and $\{e(\lambda)\}$ be as before.  Now by assumption, the function

$$z(\omega) =\max \{1/n \mid \tau(\omega)(e(1/n)(\omega)) <1 \}$$
is a.e. defined, non-zero, and finite.  It is measurable by construction, so $z$ and $z^{-1}$ represent elements of the extended center.  Now write $\varphi = (\tau_z)_{z^{-1}h}.$  Let $f$ be the spectral projection of $z^{-1}h$ for [0,1].  We have $f(\omega) = e(z(\omega))(\omega)$, so $\tau(\omega) (f(\omega)) <1$.  Then $\Phi(f) <1$, and $\tau_z(f) = \tilde{\tau} (z \Phi(f)) < \tilde{\tau}(1) <\infty.$  Since $f$ was a spectral projection of $z^{-1}h$, we conclude that $(z^{-1}h)^{-1}$ is $\tau_z$-measurable.

\end{proof}

\begin{proposition}
Let $\mathcal{M}$ be a factor with left module ${}_\mathcal{M} \mathfrak{K}$ and right module $\mathfrak{H}_\mathcal{M}$.  The map

\begin{equation}
\mathfrak{D}(\mathfrak{H},\varphi) \times \mathfrak{K} \to \mathfrak{H} \otimes_\varphi \mathfrak{K}: \quad (\xi,\eta) \mapsto \xi \otimes_\varphi \eta
\end{equation}
is preclosed only under the same conditions as in Theorem $\ref{T:preclosed}$; i.e. $\mathcal{M} = (\mathcal{M},\tau)$ is semifinite and $h^{-1}$ is $\tau$-measurable, where $\varphi = \tau_h$.
\end{proposition}

\begin{proof}
If $\mathcal{M}$ is type III, any separable $L^2$ left or right module is isomorphic to $L^2(\mathcal{M})$, so multiplication is not preclosed.

Now assume $h^{-1}$ is $\tau$-measurable.  When

$$\left(
\begin{smallmatrix}
x_1^k \varphi^{1/2} \\
x_2^k \varphi^{1/2} \\
\vdots
\end{smallmatrix}
\right)
\overset{L^2}{\underset{k \to \infty}{\longrightarrow}}
\left(
\begin{smallmatrix}
x_1 \varphi^{1/2} \\
x_2 \varphi^{1/2} \\
\vdots
\end{smallmatrix}
\right), \; 
\left(
\begin{smallmatrix}
\eta_1^k & \eta_2^k & \cdots
\end{smallmatrix}
\right) \overset{L^2}{\underset{k \to \infty}{\longrightarrow}}
\left(
\begin{smallmatrix}
\eta_1 & \eta_2 & \cdots
\end{smallmatrix}
\right), \; (x_i^k \eta_j^k) \overset{L^2}{\underset{k \to \infty}{\longrightarrow}} (\zeta_{ij}),$$
we also have $L^2$ convergence in each coordinate.  By Theorem \ref{T:preclosed}, $\zeta_{ij} = x_i \eta_j$.

When $h^{-1}$ is not $\tau$-measurable, $\mathcal{M}$ must be $\text{I}_\infty$ or $\text{II}_\infty$.  In this case $\mathcal{M} \simeq M_\infty(\mathcal{M})$, and we do not need row and column matrices: $\mathfrak{H} \simeq q_1 \, L^2(\mathcal{M})$ and $\mathfrak{K} \simeq L^2(\mathcal{M}) \, q_2$ for appropriate projections $q_1$, $q_2$.  Fix equivalent finite projections $f_1 \le q_1$, $f_2 \le q_2$ with $v^*v = f_1$, $v v^* = f_2$.  By assumption $e(1/{n^3})$ is infinite for all $n$; let $g_n$ be a subprojection equivalent to the $f_i$ with $v_{in}^* v_{in} = f_i$, $v_{in} v_{in}^* = g_n$.  Then

$$(nv_{1n}^*\varphi^{1/2}, (1/n) v_{2n}) \mapsto v^*,$$

$$||nv_{1n}^*\varphi^{1/2}||^2 = n^2 \tau(g_n h) \le (1/n) \tau(f_1) \to 0,$$

$$||(1/n) v_{2n}||^2 = (1/{n^2}) \tau(f_2) \to 0,$$
and the map is not preclosed.

\end{proof}

\end{document}